\documentclass[reqno,10pt]{amsart}

\usepackage[foot]{amsaddr}
\usepackage{graphicx,enumerate,nicefrac,bm}
\usepackage[dvipsnames]{xcolor}
\usepackage{cite}
\usepackage{amsmath, amssymb}
\usepackage{tikz}\usetikzlibrary{matrix}
\usepackage{a4wide}
\usepackage{subfig}
\usepackage{changes}
\usepackage{mathtools}
\usepackage{booktabs}
\usepackage{algorithm,algpseudocode,tabularx}
\usepackage{comment}
\usepackage{stmaryrd}

\makeatletter
\newcommand{\multiline}[1]{%
  \begin{tabularx}{\dimexpr\linewidth-\ALG@thistlm}[t]{@{}X@{}}
    #1
  \end{tabularx}
}
\makeatother

\renewcommand{\P}{\mathsf{P}}

\newcommand{\J}{\mathrm{J}}
\renewcommand{\H}{\mathsf{H}}

\newcommand{\F}{\mathsf{F}}

\newcommand{\dx}{\,\mathsf{d}\bm{x}}
\newcommand{\ds}{\,\mathsf{d}s}

\newcommand{\Rnum}[1]{\mathrm{\uppercase\expandafter{\romannumeral #1\relax}}}

\newcommand{\norm}[1]{\left\|#1\right\|}
\newcommand{\nnn}[1]{{\left\vert\kern-0.25ex\left\vert\kern-0.25ex\left\vert #1 
    \right\vert\kern-0.25ex\right\vert\kern-0.25ex\right\vert}}

\newcommand{\dprod}[1]{\langle #1\rangle}

\newcommand{\lf}{L_\F}
\newcommand{\ch}{C_\H}

\DeclareMathOperator*{\argmin}{arg\,min}

\newtheorem{theorem}{Theorem}[section]

\newtheorem{proposition}[theorem]{Proposition}

\theoremstyle{definition}

\newtheorem{remark}[theorem]{Remark}

\title[Link between the SD method and FP iterations]{A link between the steepest descent method and fixed-point iterations}

\author[P.~Heid]{Pascal Heid}
\email{pascal.heid@maths.ox.ac.uk}

\address{Mathematical Institute, University of Oxford, Woodstock Road, Oxford OX2 6GG, UK}

\thanks{The author acknowledges the financial support of the Swiss National Science Foundation (SNF), Project No. P2BEP2\underline{\space}191760.}



\begin{document}

\begin{abstract}
We will make a link between the steepest descent method for an unconstrained minimisation problem and fixed-point iterations for its Euler--Lagrange equation. In this context, we shall rediscover the preconditioned conjugate gradient method for the discretised problem. The benefit of the link between the two methods will be illustrated by a numerical experiment.
\end{abstract}

\maketitle

\section{Introduction}
Throughout this work, let $X$ be a Hilbert space endowed with an inner-product denoted by $(\cdot,\cdot)_X$ and induced norm $\norm{\cdot}_X$. Furthermore, we consider a functional $\H:X \to \mathbb{R}$ and are interested in the optimisation problem
\begin{align} \label{eq:minproblem}
\argmin_{u \in X} \H(u).
\end{align}  
In general, there may exist several (local) minimisers or possibly none at all. In this work, we shall make the following assumptions on the functional $\H$:
\begin{enumerate}[(H1)]
\item $\H$ is Gateaux-differentiable;
\item $\H$ is strictly convex;
\item $\H$ is weakly coercive, i.e., $\H(u) \to \infty$ as $\norm{u}_X \to \infty$.
\end{enumerate} 
Those assumptions (H1)--(H3) imply that $\H$ has a unique minimiser $u^\star \in X$; see, e.g.,~\cite[Thm.~25.E]{Zeidler:90}. A well-known  procedure to approximate a minimiser of such a functional $\H$ is the \emph{steepest descent method}, introduced by Augustin Cauchy in the work~\cite{CauchySD}. The main idea of this method is very intuitive: at each iteration step, we move in direction of the steepest descent. In particular, if $u^n \in X$ is a given iterate, then we set 
\begin{align} \label{eq:steepestdescentupdate}
u^{n+1}:=u^n - \delta^n \nabla \H(u^n),
\end{align}  
where $\nabla \H(u^n)$, to be specified in Section~\ref{sec:link}, is the gradient of $\H$ at $u^n$ and $\delta^n >0$ is an appropriate step-size such that $\H(u^{n+1}) \leq \H(u^n)$. The optimal choice of the step-size is given by
\begin{align*}
\delta^n:=\argmin_{t \geq 0}  \H(u^n-\delta^n \nabla \H(u^n)),
\end{align*} 
which requires the solution of a one-dimensional optimisation problem. In practice, we may often only approximate the optimal step-size. More comments on that issue will be provided in Remark~\ref{rem:pncg} below. \\

It is well-known, see, e.g.,~\cite[Thm.~25.F.]{Zeidler:90}, that under the assumptions (H1)--(H3) the unconstrained minimisation problem~\eqref{eq:minproblem} is equivalent to the operator equation
\begin{align} \label{eq:operatoreq}
\text{find} \ u \in X \ \text{such that} \qquad \F(u)=0 \qquad \text{in} \ X^\star,
\end{align}
where $\F:=\H':X \to X^\star$ is the Gateaux-derivative of the functional $\H$ and $X^\star$ denotes the dual space of $X$; i.e., the set of all continuous linear functionals from $X$ to $\mathbb{R}$. There exists a wide variety of fixed-point iterations for the numerical solution of the problem~\eqref{eq:operatoreq} and, as was shown in~\cite{HeidWihler:19v2, HeidWihler2:19v1}, in many cases they can be interpreted as an iterative local linearisation, which can be obtained by applying a suitable preconditioning operator. In particular, for any given $u \in X$, let $\P[u]:X \to X^\star$ be a linear and invertible operator. Then, the operator equation~\eqref{eq:operatoreq} is equivalent to the fixed-point equation
\begin{align} \label{eq:fpeq}
\text{find} \ u \in X \ \text{such that} \qquad u=u-\P[u]^{-1}\F(u).
\end{align} 
This, in turn, gives rise to the fixed-point iteration
\begin{align} \label{eq:fpiteration}
u^{n+1}:=u^n-\P[u^n]^{-1}\F(u^n),
\end{align} 
where $u^0 \in X$ is an initial guess. In practice, we seldom invert $\P[u^n]$ (in the finite-dimensional setting), but rather solve the linear problem
\begin{align} \label{eq:fplinear}
\text{find} \ u^{n+1} \in X \ \text{such that} \qquad \P[u^n]u^{n+1}=\P[u^n]u^n-\F(u^n) \qquad \text{in} \ X^\star
\end{align}
by applying an iterative linear solver. However, in this work, we will neglect the algebraic error for simplicity; i.e., we assume that the linear equation~\eqref{eq:fplinear} is solved exactly. Moreover, problem~\eqref{eq:fplinear} can be stated equivalently as 
\begin{align} \label{eq:fplinearalt}
\text{find} \ u^{n+1} \in X \ \text{such that} \qquad a(u^n;u^{n+1},v)=\dprod{\ell(u^n),v} \qquad \text{for all} \ v \in X,
\end{align} 
where 
\begin{align*}
a(u;v,w):=\dprod{\P[u]v,w}, \qquad u,v,w \in X
\end{align*}
and 
\begin{align*}
\dprod{\ell(u),v}:=\dprod{\P[u]u-\F(u),v}, \qquad u,v \in X;
\end{align*}
here, $\dprod{\cdot,\cdot}$ denotes the duality pairing in $X^\star \times X$. We note that some prominent iteration schemes, such as the Zarantonello, Ka\v{c}anov, and Newton methods, can be cast into this unified framework, we refer to~\cite{HeidWihler:19v2}. In the following we assume that, for any fixed $u \in X$, the bilinear form $a(u;\cdot,\cdot):X \times X \to \mathbb{R}$ is 
\begin{enumerate}[(A1)]
\item[(A1)] uniformly coercive, i.e., there exists a constant $\alpha>0$ such that
\begin{align*}
a(u;v,v) \geq \alpha \norm{v}_X^2 \qquad \text{for all} \ u,v \in X;
\end{align*}
\item[(A2)] uniformly bounded, i.e., there exists a constant $\beta>0$ such that
\begin{align*}
a(u;v,w) \leq \beta \norm{v}_X \norm{w}_X \qquad \text{for all} \ u,v,w \in X;
\end{align*}
\item[(A3)] symmetric, i.e., $a(u;v,w)=a(u;w,v)$ for all $u,v,w \in X$.
\end{enumerate}
We note that those assumptions imply, thanks to the Lax--Milgram theorem, that the operator $\P[u]$ is invertible for any $u \in X$ and that the equation~\eqref{eq:fplinearalt} has a unique solution for each $n=0,1,2,\dotsc$. We further assume that the operator $\F:X \to X^\star$ is 
\begin{enumerate}[(F1)]
\item strongly monotone, i.e., there exists a constant $\nu>0$ such that
\begin{align*}
\dprod{\F(u)-\F(v),u-v} \geq \nu \norm{u-v}^2_X \qquad \text{for all} \ u,v \in X;
\end{align*}
\item Lipschitz continuous, i.e., there exists a constant $\lf>0$ such that
\begin{align*}
\dprod{\F(u)-\F(v),w} \leq \lf \norm{u-v}_X \norm{w}_X \qquad \text{for all} \ u,v,w \in X.
\end{align*}
\end{enumerate}
Under those assumptions, the theory on strongly monotone operator equations yields that equation~\eqref{eq:operatoreq} has a unique solution $u^\star \in X$; see, e.g.,~\cite[\S 25.4]{Zeidler:90}. We recall that this solution is as well the unique minimiser of $\H$ in $X$. We further note that the strong monotonicity (F1) of the operator $\F$ implies the strict convexity (H2) of its potential $\H$.\\

Under certain assumptions it can be shown that the potential $\H$ decreases along the sequence $\{u^n\}$ generated by the unified iteration scheme~\eqref{eq:fplinearalt} in the sense that there exists a constant $\ch>0$ such that
\begin{align} \label{eq:energydecay}
\H(u^n)-\H(u^{n+1}) \geq \ch \norm{u^n-u^{n+1}}_X^2 \qquad \text{for all} \ n=0,1,2,\dotsc;
\end{align}
see~\cite[\S 2.1]{HeidWihler2:19v1} for a general discussion of the property~\eqref{eq:energydecay} and~\cite[\S 2.4]{HeidPraetoriusWihler:21} for the required assumptions guaranteeing~\eqref{eq:energydecay} for the Zarantonello, Ka\v{c}anov, and Newton methods. In particular, given the monotonicity property~\eqref{eq:energydecay}, the update $-\P[u^n]^{-1}\F(u^n)$ from the fixed-point iteration~\eqref{eq:fpiteration} can be considered as a descent direction of the potential $\H$ at the given iterate $u^n \in X$. This indicates that there might be a link between the steepest descent method~\eqref{eq:steepestdescentupdate} for the optimisation problem~\eqref{eq:minproblem} and the fixed-point iteration~\eqref{eq:fpiteration} for the solution of its Euler--Lagrange equation~\eqref{eq:operatoreq}. The purpose of this work is to make this connection visible.

\subsection*{Outline} In Section~\ref{sec:link} we will make a link between the steepest descent method~\eqref{eq:steepestdescentupdate} and the unified iteration scheme~\eqref{eq:fpiteration}. Subsequently, in Section~\ref{sec:PNCG}, we consider the nonlinear conjugate gradient method with regard to our insights from Section~\ref{sec:link}. A numerical experiment is performed in Section~\ref{sec:experiments} and finally we will round off our work with some conclusions in Section~\ref{sec:con}.  

\section{Link between the steepest descent method and the unified iteration scheme} \label{sec:link}

We will now make the link between the steepest descent method~\eqref{eq:steepestdescentupdate} and the unified iteration scheme~\eqref{eq:fpiteration} visible. For that purpose, let us first recall the steepest descent method~\eqref{eq:steepestdescentupdate}, which involves the gradient $\nabla \H(u^n) \in X$ of $\H$ at $u^n \in X$. By definition it holds that, for fixed $u \in X$,
\begin{align} \label{eq:steepestd}
\dprod{\F(u),v}=\dprod{\H'(u),v}=:(\nabla \H(u),v)_X \qquad \text{for all} \ v \in X.
\end{align}
In particular, the gradient depends on the considered inner-product, which shall be indicated by a subscript in the following; i.e., we write $\nabla_X \H(u)$ for the gradient of $\H$ at $u \in X$ with respect to the inner-product $(\cdot,\cdot)_X$. Moreover, if we denote by $\J_X:X \to X^\star$ the Riesz isometry with respect to the inner-product $(\cdot,\cdot)_X$ on $X$, then we have that $\nabla_X \H(u)=\J_X^{-1} \H'(u)=\J_X^{-1} \F(u)$, cf.~\eqref{eq:steepestd}. In turn, the steepest descent method reads as
\begin{align} \label{eq:steepestdescent2}
u^{n+1}=u^n-\delta^n\J_X^{-1}\F(u^n).
\end{align}
Consequently, the steepest descent method~\eqref{eq:steepestdescent2} coincides with the fixed-point iteration~\eqref{eq:fpiteration} for the preconditioning operator $\P[u]={\delta(u)}^{-1}\J_X$, $u \in X$, where the damping function satisfies $\delta(u^n)=\delta^n$ for $n=0,1,2,\dotsc$. We note that this specific choice of the preconditioning operator gives rise to the Zarantonello iteration; see the original work~\cite{Zarantonello:60}, or the monographs~\cite[\S3.3]{Necas:86} and~\cite[\S25.4]{Zeidler:90}. Moreover, given the assumptions (F1)--(F2), the Zarantonello iteration generates a sequence converging to the unique solution $u^\star \in X$ of~\eqref{eq:operatoreq} for a suitable choice of the damping function $\delta:X \to \mathbb{R}_{>0}$; see, e.g.,~the proof of \cite[Thm.~25.B]{Zeidler:90}. 

If $a:X \times X \to \mathbb{R}$ is a symmetric, coercive, and bounded bilinear form on $X \times X$, then, in particular, $a(\cdot,\cdot)$ is an inner-product on $X$ whose corresponding norm $\norm{\cdot}_a$ is equivalent to the norm $\norm{\cdot}_X$; i.e., $X$ endowed with the inner-product $a(\cdot,\cdot)$ and norm $\norm{\cdot}_a$ is a Hilbert space as well. We further note that, in turn, the bilinear form $a:X \times X \to \mathbb{R}$ induces a linear and invertible operator $\P:X \to X^\star$ defined by
\begin{align} \label{eq:Afroma}
\dprod{\P u,v}:=a(u,v) \qquad \text{for all} \ u,v \in X.
\end{align}
We may then consider the gradient with respect to the inner-product $a(\cdot,\cdot)$, i.e., for given $u \in X$,
\begin{align*}
a(\nabla_a \H(u),v):=\dprod{\H'(u),v} \qquad \text{for all} \ v \in X.
\end{align*}
In view of~\eqref{eq:Afroma} we have that
\begin{align*}
\dprod{\P \nabla_a \H(u),v}=\dprod{\H'(u),v}=\dprod{\F(u),v} \qquad \text{for all} \ v \in X
\end{align*}
and therefore $\nabla_a \H(u)=\P^{-1} \F(u)$. In this case, the steepest descent method~\eqref{eq:steepestdescentupdate} coincides, up to some damping parameter, with the unified iteration scheme~\eqref{eq:fpiteration} for the preconditioner from~\eqref{eq:Afroma}.

Finally, similarly as was done in~\cite{HenningPeterseim:18} in the context of Sobolev gradient flows for the Gross--Pitaevskii equation, we may consider an inner-product that changes with the iteration. For fixed $u \in X$, let $a_u=a(u;\cdot,\cdot):X \times X \to \mathbb{R}$ be a symmetric, uniformly coercive and bounded bilinear form, cf.~(A1)--(A3). Consequently, for any $u \in X$, the operator $\P[u]: X \to X^\star$ defined by
\begin{align} \label{eq:Aufromau}
\dprod{\P[u]v,w}:=a_u(v,w)=a(u;v,w), \qquad v,w \in X,
\end{align} 
is linear and invertible. Then, we can define the gradient of $\H$ at a given element $u \in X$ with respect to the inner-product $a_u(\cdot,\cdot)$ by
\begin{align} \label{eq:changinggradient}
a_u(\nabla_{a_u} \H(u),v):=\dprod{\H'(u),v}=\dprod{\F(u),v} \qquad \text{for all} \ v \in X;
\end{align}
i.e., we have that $\nabla_{a_u} \H(u)=\P[u]^{-1}\F(u)$. In turn, the steepest descent method is given by
\begin{align*}
u^{n+1}=u^n-\delta(u^n)\P[u^n]^{-1}\F(u^n),
\end{align*}
which, for $\delta(u^n)=1$, $n=0,1,2,\dotsc$, matches our unified iteration scheme~\eqref{eq:fpiteration}. In particular, we have shown the following result.

\begin{proposition} 
Let $\P[u]:X \to X^\star$, for $u \in X$, be a linear and invertible operator which induces a bilinear form $a_u=a(u;\cdot,\cdot)$ that satisfies {\rm (A1)--(A3)} (or vice versa), cf.~\eqref{eq:Aufromau}. Then, the unified iteration scheme~\eqref{eq:fpiteration} with preconditioner $\P[u]$ coincides with the steepest descent method~\eqref{eq:steepestdescentupdate} where the gradient is taken with respect to the (changing) inner-product $a_u(\cdot,\cdot)$ and the step-sizes satisfy $\delta^n = 1$ for all $n=0,1,2,\dotsc$. 
\end{proposition}

\begin{remark}
There may result some advantages from this connection between the steepest descent method and the unified fixed-point iteration. 
\begin{enumerate}[(1)]
\item The step-size function of the (modified) steepest descent method in this context is simply given by $\delta \equiv 1$ and thus we do not need to employ, e.g., a line search or trusted region method to determine $\delta^n$, $n=0,1,2,\dotsc$. We note that the preconditioning operator $\P[u]$ may implicitly include a damping parameter $\delta(u)$; however, in many cases, this damping parameter can be prescribed or can easily be chosen adaptively in such a way that the decay property~\eqref{eq:energydecay} is satisfied in each iteration step. 
\item There is a large body of literature focusing on the convergence of fixed-point iterations. By the identification of the unified iteration scheme~\eqref{eq:fpiteration} and the steepest descent method~\eqref{eq:steepestdescentupdate} (with constant step-size function $\delta \equiv 1$), those results also apply to the latter.
\item On the other hand, the steepest descent method serves as basis of the superior (nonlinear) conjugate gradient method. Hence, it might be sensible to consider the nonlinear conjugate gradient method in the case that the gradient is taken with respect to an inner-product induced by a preconditioning operator $\P[u]:X \to X^\star$ (motivated by a fixed-point iteration), cf.~\eqref{eq:changinggradient}. Indeed, we will show in the next section that this gives rise to the known preconditioned nonlinear conjugate gradient (PNCG) method. 
\end{enumerate}
\end{remark}

\section{Nonlinear conjugate gradient method for general inner-products} \label{sec:PNCG}

For the purpose of examining the nonlinear conjugate gradient method in the context of dynamic inner-products we will consider a model problem, which will now be introduced. 

\subsection{Model problem} 

As our model problem, we will consider the following quasilinear second-order elliptic partial differential equation:  
\begin{align} \label{eq:modelproblem}
\text{find} \ u \in X \ \text{such that} \qquad -\nabla \cdot \{\mu(|\nabla u|^2)\nabla u\}-g =0  \qquad \text{in } X^\star,
\end{align}
i.e., we set 
\begin{align} \label{eq:modelF}
\F(u):=-\nabla \cdot \{\mu(|\nabla u|^2)\nabla u\}-g
\end{align}
in~\eqref{eq:operatoreq}. Here, $X:=H_0^1(\Omega)$ is the Sobolev space of $H^1$-functions with zero trace along the boundary $\partial \Omega$, where $\Omega \subset \mathbb{R}^d$, $d \in \{2,3\}$, is an open, bounded, and polygonal domain. For $u,v \in X$, the inner-product and norm on $X$ are defined by $(u,v)_X:=(\nabla u,\nabla v)_{L^2(\Omega)}$ and $\norm{u}_X:=\norm{\nabla u}_{L^2(\Omega)}$, respectively. Furthermore, $g \in L^2(\Omega)$, considered as an element in the dual space $H^{-1}(\Omega):=H_0^1(\Omega)^\star$, is given and the diffusion coefficient $\mu \in C^1([0,\infty))$ satisfies the monotonicity condition
\begin{align} \label{eq:muass}
m_\mu (t-s) \leq \mu(t^2)t-\mu(s^2)s \leq M_\mu (t-s), \qquad t \geq s \geq 0,
\end{align}
for some constants $M_\mu \geq m_\mu >0$. Under those assumptions, the nonlinear operator $\F:X \to X^\star$ from~\eqref{eq:modelF} satisfies the conditions (F1) and (F2) with $\nu=m_\mu$ and $\lf=3 M_\mu$; see, e.g.,~\cite[Prop.~25.26]{Zeidler:90}. We note the weak form of our model problem~\eqref{eq:modelproblem}:
\begin{align} \label{eq:modelproblemweak}
\text{find} \ u \in X \ \text{such that} \qquad \int_\Omega \mu(|\nabla u|^2)\nabla u \cdot \nabla v \dx=\int_\Omega gv \dx \qquad \text{for all} \ v \in X.
\end{align}
It is straightforward to verify that $\F$ is a potential operator with the potential given by 
\begin{align*}
\H(u):=\int_\Omega \psi(|\nabla u|^2) \dx -\int_\Omega gu \dx, \qquad u \in X,
\end{align*}
where $\psi(s)=\nicefrac{1}{2} \int_0^s \mu(s) \ds$ for $s \geq 0$. As $\F$ is strongly monotone it immediately follows that $\H$ is strictly convex. Furthermore, the Cauchy--Schwarz inequality, the Poincar\'{e}--Friedrich inequality (with constant denoted by $C_P$), and the assumption~\eqref{eq:muass} imply that
\begin{align*}
\H(u) \geq \frac{m_\mu}{2} \norm{u}^2_{X}-C_P \norm{g}_{L^2(\Omega)}\norm{u}_X,
\end{align*}
thus $\H$ is weakly coercive. In particular, (H1)--(H3) are satisfied. If we further assume that $\mu$ is monotonically decreasing, i.e., $\mu'(t) \leq 0$ for all $t \geq 0$, then the Zarantonello iteration, the Ka\v{c}anov scheme, and the damped Newton method all satisfy --- for suitable damping functions ---
assumptions (A1)--(A3) as well as~\eqref{eq:energydecay}; see~\cite{HeidWihler:19v2,HeidWihler2:19v1,HeidPraetoriusWihler:21}. In particular, the following three methods generate a sequence converging to the unique solution of~\eqref{eq:operatoreq} and, equivalently, of~\eqref{eq:minproblem}:
\begin{enumerate}[(i)]
\item \emph{Zarantonello (or Picard) iteration,} for $\delta_Z \in (0,\nicefrac{2}{3 M_\mu})$:
\begin{align} \label{eq:zarantonello}
u^{n+1}=u^n-\delta_Z \J_X^{-1}\F(u^n) \qquad \text{for all} \ n=0,1,2,\dotsc,
\end{align}
where $\J_X:X \to X^\star$ denotes, as before, the Riesz isometry with respect to the inner-product $(\cdot,\cdot)_X$ on $X$;
\item \emph{Ka\v{c}anov iteration:}
\begin{align*} 
u^{n+1}=u^{n}-\P[u^n]^{-1}\F(u^n)\qquad \text{for all} \ n=0,1,2,\dotsc,
\end{align*}
where $\dprod{\P[u]v,w}:=\int_\Omega \mu(|\nabla u|^2) \nabla v \cdot \nabla w\dx$ for $u,v,w \in X$, or, equivalently,
\begin{equation*} 
- \nabla \cdot \big\{\mu(\left|\nabla u^n\right|^2) \nabla{u^{n+1}}\big\}=g \qquad \text{for all} \ n=0,1,2,\dotsc;
\end{equation*}
\item \emph{Newton iteration,} for a damping parameter $0 < \delta_{\mathrm{min}} \leq \delta_N(u^n) \leq \delta_{\mathrm{max}}<\nicefrac{2 m_\mu}{3 M_\mu}$:
\begin{align} \label{eq:newton}
u^{n+1}=u^{n}-\delta_N(u^n)  \F'(u^{n})^{-1}\F(u^{n})\qquad \text{for all} \ n=0,1,2,\dotsc;
\end{align}
here, for $u\in X$, the Gateaux-derivative $\F'(u)$ of $\F$ is given through
\[
\dprod{\F'(u)v,w}=\int_{\Omega} 2 \mu'(|\nabla u|^2)(\nabla u \cdot \nabla v)(\nabla u \cdot \nabla w) \dx + \int_{\Omega} \mu(|\nabla u|^2)\nabla v \cdot \nabla w \dx, \quad v,w \in X.
\]
\end{enumerate}

\subsection{Discretisation of the model problem} \label{sec:discretemodel} Since $X=H_0^1(\Omega)$ is an infinite-dimensional space, we cannot compute the sequence generated by any of the iteration schemes presented before. In order to cast them into a computational framework, we will consider the discretisation by the conforming $\mathbb{P}_1$-finite element method. In particular, let $\mathcal{T}$ be a triangulation of $\Omega$ and the corresponding $\mathbb{P}_1$-finite element space is given by
\[X_h:=\{u \in H_0^1(\Omega): u|_K \in \mathbb{P}_1(K) \ \forall K \in \mathcal{T}\},\]
where $\mathbb{P}_1(K)$ denotes the set of all affine functions on $K$. Then, the discretisation of the weak problem~\eqref{eq:modelproblemweak} reads as
\begin{align} \label{eq:discreteweak}
\text{find} \ u_h \in X_h \ \text{such that} \qquad \int_\Omega \mu(|\nabla u_h|^2)\nabla u_h \cdot \nabla v \dx =\int_\Omega g v \dx \qquad \text{for all} \ v \in X_h.
\end{align}  
Furthermore, upon defining 
\begin{align*} 
\mathsf{B}_h(u;v,w):=\int_\Omega \mu(|\nabla u|^2) \nabla v \cdot \nabla w \dx, \qquad u,v,w \in X_h,
\end{align*}
and
\begin{align*} 
\dprod{\ell_h,v}:=\int_\Omega gv\dx, \qquad v \in X_h,
\end{align*}
the discrete weak problem~\eqref{eq:discreteweak} can be stated equivalently as follows:
\begin{align} \label{eq:discreteweak2}
\text{find} \ u_h \in X_h \ \text{such that} \qquad \mathsf{B}_h(u_h;u_h,v) =\dprod{\ell_h,v} \qquad \text{for all} \ v \in X_h.
\end{align}
We emphasise that, for any $u \in X_h$, $\mathsf{B}_h(u;\cdot,\cdot):X_h \times X_h \to \mathbb{R}$ is a symmetric, uniformly coercive and bounded bilinear form. In particular, we have that
\begin{align*}
\mathsf{B}_h(u;v,v) \geq m_\mu \norm{v}_X^2 \qquad \text{for all} \ u,v \in X_h
\end{align*}
and
\begin{align*}
\mathsf{B}_h(u;v,w) \leq M_\mu \norm{v}_X \norm{w}_X \qquad \text{for all} \ u,v,w \in X_h.
\end{align*}
Consequently, thanks to the Lax--Milgram theorem, \eqref{eq:discreteweak2} has a unique solution. Moreover, if we define $\mathsf{F}_h: X_h \to X_h^\star$ by
\begin{align} \label{eq:Fmatrix}
\mathsf{F}_h(u):=\mathsf{B}_h(u;u,\cdot)-\ell_h, \qquad u \in X_h,
\end{align}
then we can state~\eqref{eq:discreteweak2} in form of an operator equation:
\begin{align*} 
\text{find} \ u_h \in X_h \ \text{such that} \qquad \mathsf{F}_h(u_h)=0 \qquad \text{in} \ X_h^\star.
\end{align*}

Now let $\{\xi_i\}_{i=1}^{m_h}$ be the nodal basis of $X_h$, where $m_h \in \mathbb{N}$ is the number of degrees of freedom in $X_h$. Consequently, each element $u \in X_h$ can be written in a unique way as a linear combination of those basis vectors; i.e., $u=\sum_{i=1}^{m_h} c_i \xi_i$, where $c_i \in \mathbb{R}$, for $i \in \{1,\dotsc,m_h\}$, are the coefficients of $u$ with respect to the basis $\{\xi_i\}_{i=1}^{m_h}$. Then, the corresponding linear mapping $\Psi:\mathbb{R}^{m_h} \to X_h$ defined by $\Psi(\mathbf{u}):=\sum_{i=1}^{m_h} c_i \xi_i$, where $\mathbf{u}=(c_1,\dotsc,c_{m_h})^{T}$, is one-to-one. By invoking this isomorphism we may consider the discrete weak equation~\eqref{eq:discreteweak2} as a problem in $\mathbb{R}^{m_h}$: 
\begin{align} \label{eq:rm1}
\text{find} \ \mathbf{u}_h \in \mathbb{R}^{m_h} \ \text{such that} \qquad \mathsf{B}_h(\Psi(\mathbf{u}_h);\Psi(\mathbf{u}_h),\Psi(\mathbf{v})) =\dprod{\ell_h,\Psi(\mathbf{v})} \qquad \text{for all} \ \mathbf{v} \in \mathbb{R}^{m_h}.
\end{align}
We note that, for any $\mathbf{u} \in \mathbb{R}^{m_h}$, $\mathsf{B}_h(\Psi(\mathbf{u}),\Psi(\cdot),\Psi(\cdot)):\mathbb{R}^{m_h} \times \mathbb{R}^{m_h} \to \mathbb{R}$ is a symmetric and coercive bilinear form on $\mathbb{R}^{m_h} \times \mathbb{R}^{m_h}$ and thus can be represented by a symmetric positive definite matrix $\mathbf{A}_h^{\mu}(\mathbf{u}) \in \mathbb{R}^{m_h\times m_h}$ (which depends on $\mathbf{u} \in \mathbb{R}^{m_h}$). Likewise, $\dprod{\ell_h,\Psi(\cdot)}:\mathbb{R}^{m_h} \to \mathbb{R}$ is a linear form and hence can be identified with a vector $\mathbf{b}_h \in \mathbb{R}^{m_h}$. Consequently, problem~\eqref{eq:rm1} can be restated as: 
\begin{align*} 
\text{find} \ \mathbf{u}_h \in \mathbb{R}^{m_h} \ \text{such that} \qquad \mathbf{A}_h^{\mu}(\mathbf{u}_h) \cdot \mathbf{u}_h = \mathbf{b}_h;
\end{align*}
here and in the following, in the context of matrices and vectors in $\mathbb{R}^{m_h \times m_h}$ and $\mathbb{R}^{m_h}$, respectively, we denote by $\cdot$ the usual matrix-vector product.

\subsection{Algebraic gradient descent method and preconditioning}
By invoking the isomorphism $\Psi: \mathbb{R}^{m_h} \to X_h$, the operator $\F_h:X_h \to X_h^\star$ from~\eqref{eq:Fmatrix} can be considered as an operator $\mathbf{F}_h:\mathbb{R}^{m_h} \to \mathbb{R}^{m_h}$ given by
\begin{align*}
\mathbf{F}_h(\mathbf{u}):=\mathbf{A}_h^{\mu}(\mathbf{u}) \cdot \mathbf{u}- \mathbf{b}_h, \qquad \mathbf{u} \in \mathbb{R}^{m_h}.
\end{align*}
In particular, this is the algebraic gradient with respect to the Euclidean inner-product on $\mathbb{R}^{m_h}$ and the corresponding gradient descent method reads as 
\[\mathbf{u}^{n+1}=\mathbf{u}^n-\delta(\mathbf{u}^n)\mathbf{F}_h(\mathbf{u}), \qquad n=0,1,2,\dotsc,\] 
where $\mathbf{u}^0 \in \mathbb{R}^{m_h}$ is an initial guess and $\delta(\mathbf{u}^n)>0$, for $n=0,1,2,\dotsc$, are a suitable step-sizes. We point out that the gradient chosen that way is completely detached from the original partial differential equation (arising as the mathematical model of, e.g., a physical problem). Therefore, we should rather consider the discrete (vector) version of the generalised gradient from~\eqref{eq:changinggradient}. In particular, for given $\mathbf{u} \in \mathbb{R}^{m_h}$, let $\mathbf{\nabla_{a_u} H(u)} \in \mathbb{R}^{m_h}$ be such that
\begin{align*}
a_{\Psi(\mathbf{u})}(\Psi(\mathbf{\nabla_{a_u} H(u)}),\Psi(\mathbf{v}))=\dprod{\F_h(\Psi(\mathbf{u})),\Psi(\mathbf{v})}=\mathbf{v}^{T} \cdot \mathbf{F}_h(\mathbf{u}) \qquad \text{for all} \ \mathbf{v} \in \mathbb{R}^{m_h}.
\end{align*}
Since $a_{\Psi(\mathbf{u})}(\Psi(\cdot),\Psi(\cdot)):\mathbb{R}^{m_h} \times \mathbb{R}^{m_h} \to \mathbb{R}$ is a symmetric and coercive bilinear form, it can be represented by a symmetric positive definite matrix $\mathbf{P}_h(\mathbf{u})\in \mathbb{R}^{m_h \times m_h}$. Hence, we have that 
\begin{align*}
\mathbf{\nabla_{a_u} H(u)}=\mathbf{P}_h(\mathbf{u})^{-1} \cdot \mathbf{F}_h(\mathbf{u})
\end{align*} 
and, in turn,
\begin{align*}
\mathbf{u}^{n+1}=\mathbf{u}^n-\delta(\mathbf{u}^n)\mathbf{P}_h(\mathbf{u})^{-1} \cdot \mathbf{F}_h(\mathbf{u}), \qquad n=0,1,2,\dotsc.
\end{align*}
In particular, if $\mathbf{P}_h(\mathbf{u})=\mathbf{P}_h$ is independent of $\mathbf{u} \in \mathbf{R}^{m_h}$, then this procedure coincides with the (algebraic) gradient descent method for the preconditioned problem
\begin{align*}
\mathbf{P}_h^{-1} \cdot \mathbf{A}_h^{\mu}(\mathbf{u}_h) \cdot \mathbf{u}_h = \mathbf{P}_h^{-1} \cdot \mathbf{b}_h.
\end{align*}
We note that this is not a new insight, but is already known in the literature: for linear problems, this and many more observations are well presented in~\cite{MalekStrakos:15}.

\subsection{Preconditioned nonlinear conjugate gradient method}
As we have seen before, at least for our model problem, the gradient descent method with respect to an inner-product $a(\cdot,\cdot)$ (induced by a linear and invertible operator $\P:X \to X^\star$) simply leads to the preconditioned (algebraic) gradient descent method. Consequently, if we want to derive the conjugate gradient method in the case that the gradient is taken with respect to some inner-product $a_u(\cdot,\cdot)$, $u \in X$, on $X$, this simply leads to the known preconditioned nonlinear conjugate gradient method, see Algorithm~\ref{alg:pcng}. More details about (the derivation of) this method can be found, e.g., in the book~\cite{Pytlak:09}, but we also refer to the article~\cite{Caliciotti:17}.  

\begin{algorithm}
\caption{Preconditioned nonlinear conjugate gradient method}
\label{alg:pcng}
\begin{algorithmic}[1]
\State Input initial guess $\mathbf{u}^0 \in \mathbb{R}^{m_h}$.
\State Compute $\mathbf{d}^0=-\mathbf{P}_h(\mathbf{u}^0)^{-1} \cdot \mathbf{F}_h(\mathbf{u}^0)$ and set $n=0$.
\Repeat 
\State Compute the step-length $\alpha(\mathbf{u}^n)$ and set $\mathbf{u}^{n+1}=\mathbf{u}^{n}+\alpha(\mathbf{u}^n)\mathbf{d}^n$.
\State Compute $\beta^n$ and set $\mathbf{d}^{n+1}=-\mathbf{P}_h(\mathbf{u}^{n+1})^{-1} \cdot \mathbf{F}_h(\mathbf{u}^{n+1})+\beta^n \mathbf{d}^n$.
\State Update $n \leftarrow n+1$.
\Until{stopping criterion is satisfied.}
\State \textsc{Return} approximate solution $\mathbf{u}^{n}$.
\end{algorithmic}
\end{algorithm}

\begin{remark} \label{rem:pncg}
Without going into too much details, we shall provide some comments on Algorithm~\ref{alg:pcng}.
\begin{enumerate}[(1)]
\item Ideally, the step-size $\alpha(\mathbf{u}^n)\geq 0$ is chosen such that 
\begin{align} \label{eq:optstep}
\alpha(\mathbf{u}^n)=\argmin_{\alpha \geq 0} \H(\Psi(\mathbf{u}^n+\alpha\mathbf{d}^n)).
\end{align}
In practice, however, we can, in general, only approximate this minimiser by using, e.g., a line search or a trusted region method; we refer, e.g., to~\cite[\S 3.2]{Kelley:99}. Often, especially for convergence proofs, it is required that the choice of the step-size satisfies some version of the Wolfe conditions. The standard Wolfe conditions were introduced in~\cite{Wolfe1,Wolfe2}. Later on, several modified Wolfe conditions were presented; we refer to~\cite{HagerZhang:06} and the references therein. 
\item Many different choices for the conjugate gradient update parameter $\beta^n$ have been proposed in the literature, see, e.g., the extensive survey of Hager and Zhang~\cite{HagerZhang:06} and the references therein. For the PNCG method, two of the most popular choices are the ones proposed by Fletcher and Reeves~\cite{FletcherReeves:64},
\begin{align} \label{eq:FR}
\beta_{FR}^n=\frac{\mathbf{F}_h(\mathbf{u}^{n+1})^T \cdot \mathbf{P}_h(\mathbf{u}^{n+1})^{-1}\cdot \mathbf{F}_h(\mathbf{u}^{n+1})}{\mathbf{F}_h(\mathbf{u}^{n})^T \cdot \mathbf{P}_h(\mathbf{u}^{n})^{-1} \cdot \mathbf{F}_h(\mathbf{u}^{n})}, 
\end{align} 
and by Polak and Ribi\`{e}re~\cite{PolakRibiere:69} and Polyak~\cite{Polyak:69},
\begin{align*} 
\beta_{PR}^n=\frac{[\mathbf{F}_h(\mathbf{u}^{n+1})-\mathbf{F}_h(\mathbf{u}^{n})]^T \cdot \mathbf{P}_h(\mathbf{u}^{n+1})^{-1} \cdot \mathbf{F}_h(\mathbf{u}^{n+1})}{\mathbf{F}_h(\mathbf{u}^{n})^T \cdot \mathbf{P}_h(\mathbf{u}^{n})^{-1} \cdot \mathbf{F}_h(\mathbf{u}^{n})}.  
\end{align*}
Later on, Powell proposed in the article~\cite{Powell:84} the following modified (and improved) version of the parameter $\beta^n_{PR}$:
\begin{align} \label{eq:PRP}
\beta_{PR+}^n=\max\{\beta_{PR}^n,0\}.
\end{align}
\end{enumerate}
\end{remark}

\section{Numerical experiment} \label{sec:experiments}
In this section we run a numerical experiment to examine the influence of the specific (operator) preconditioner on the performance of the nonlinear conjugate gradient method. To this end we consider our model problem~\eqref{eq:modelproblemweak}, where $\Omega:=(-1,1)^2 \setminus [0,1] \times [0,1] \subset \mathbb{R}^2$ is an L-shaped domain and the diffusion coefficient $\mu$ obeys the Carreau law; i.e., we have that
\begin{align*} 
\mu(t)=\mu_\infty+(\mu_0-\mu_\infty)(1+\lambda t)^{\nicefrac{(r-2)}{2}},
\end{align*} 
with $\mu_0>\mu_\infty>0$, $\lambda>0$, and $r \in (1,2)$. It is straightforward to verify that this choice of the diffusion coefficient satisfies~\eqref{eq:muass} with $m_\mu=\mu_\infty$ and $M_\mu=\mu_0$. Moreover, since $r \in (1,2)$, the diffusion coefficient is decreasing. Therefore, the Zarantonello, Ka\v{c}anov, and Newton methods converge for appropriate damping parameters. The source term $g \in L^2(\Omega)$ is chosen such that the unique solution of~\eqref{eq:modelproblemweak} is given by the smooth function 
\begin{align*} 
u^\star(x,y)=\sin(\pi x)\sin(\pi y),
\end{align*} 
where $(x,y) \in \mathbb{R}^2$ denote the Euclidean coordinates.
Furthermore, for the discretisation of problem~\eqref{eq:modelproblemweak}, we consider, as in Section~\ref{sec:discretemodel}, the conforming $\mathbb{P}_1$-finite element method, where the mesh $\mathcal{T}$ consists of $\mathcal{O}(10^5)$ triangles. In our experiment below, we choose the parameters $\mu_\infty=1$, $\mu_0=100$, $\lambda=2$, and (a) $r=1.4$ or (b) $r=1.05$. 
In order to approximate the corresponding solutions of the discretised problem~\eqref{eq:discreteweak} for the parameters from (a) and (b), respectively, we will apply the Ka\v{c}anov method with 1000 iteration steps. Subsequently, we will examine how many iteration steps are required by the Zarantonello, Ka\v{c}anov, and Newton methods, as well as their conjugated counterparts with update parameters from~\eqref{eq:FR} and \eqref{eq:PRP}, respectively, in order to obtain an error-tolerance of $10^{-6}$ with respect to the norm $\norm{\cdot}_X$ in $X$. In each case we choose the function $u^0 \equiv 0 \in X_h$ as our initial guess. Moreover, the one-dimensional optimisation problem from line 4 in Algorithm~\ref{alg:pcng}, cf.~\eqref{eq:optstep}, is solved by the Matlab subroutine \emph{fmincon} from the optimisation toolbox.   

In Table~\ref{table:steps} we record the number of iteration steps that were performed by our nonlinear solvers to obtain an error-tolerance of $10^{-6}$. If this accuracy was not achieved in 100 iteration steps, then the calculations were aborted, signified by '-' in the table below. The damping parameters in the Zarantonello iteration~\eqref{eq:zarantonello} were chosen to be $\delta_Z=0.03$ in (a) and $\delta_Z=0.02$ in (b), respectively, as they seemed to be close to optimal. Morever, in both cases we set $\delta_N \equiv 1$ in~\eqref{eq:newton}, i.e., we considered the classical (undamped) Newton method. We further note that neither the algebraic gradient descent nor the conjugate gradient method (without preconditioning) converged in a reasonable number of iteration steps; hence, they are not included in the table below.

\begin{table}[h] 
\begin{tabular}{l|ccc|ccc}\toprule
& \multicolumn{3}{c}{(a) $r=1.4$} & \multicolumn{3}{c}{(b) $r=1.05$}
\\ 
           & FP & $\beta_{FR}^n$  & $\beta_{PR+}^n$    & FP & $\beta_{FR}^n$  & $\beta_{PR+}^n$ \\\midrule
Zarantonello  & 61 & 15 & 15 & - & 37 & 37 \\
Ka\v{c}anov & 25 & 9 & 10 & 90 & 19 & 24 \\
Newton & 5 & 7 & 6 & 7 & 16 & 8\\\bottomrule
\end{tabular}
\caption{The required number of iteration steps for the different nonlinear solvers to obtain an error tolerance of $10^{-6}$. Here, 'FP' signifies the usual fixed-point iteration. If the prescribed accuracy was not obtained in 100 steps, then this is remarked by the sign '-' in the table.}
\label{table:steps}
\end{table}

As we can see in Table~\ref{table:steps}, the Newton method outperformed the other iteration schemes in the specific problem considered. Indeed, the classical Newton method was even slightly superior to its conjugated counterparts. In contrast, considering the Ka\v{c}anov and Zarantonello schemes, we observe that their corresponding PNCG methods require significantly less iteration schemes, at least in the given experiment. However, we should be aware that the preconditioned conjugate gradient methods require an additional solution of a one-dimensional minimisation problem, which, in general, is not for free. 

Finally, we note that for more complicated problems the domain of convergence for the Newton scheme might be quite small, and thus we have to consider other nonlinear solvers such as the Ka\v{c}anov and Zarantonello methods; see, e.g.,~\cite[\S 5.1]{HeidSuli2:21}. Hence, it is certainly worth to study those iteration schemes. Moreover, as we have seen above, their conjugated counterparts are able to accelerate the convergence (in view of the number of iteration steps), at least for the model problem considered.    

\section{Conclusion} \label{sec:con}

We showed that, up to a damping function, the fixed-point iteration obtained by a preconditioning operator coincides with the steepest descent method in the case that the gradient is taken with respect to the inner-product induced by this preconditioning operator. Moreover, in view of the corresponding discretised problem in $\mathbb{R}^{m_h}$, the operator preconditioner acts as an algebraic preconditioner and, in turn, leads to the preconditioned gradient descent method. Our numerical experiment illustrated that the choice of a problem related (operator) preconditioner may significantly improve the convergence of the nonlinear conjugate gradient method.  
 
\bibliographystyle{amsplain}
\bibliography{references}
\end{document}